\documentclass[10pt]{wlscirep}
\usepackage[utf8]{inputenc}
\usepackage[T1]{fontenc}
\usepackage{todonotes}
\usepackage{graphicx}
\usepackage{xcolor}
\usepackage{wrapfig}
\usepackage{soul}
\usepackage[toc,page]{appendix}

\def\R{\mathbb{R}}
\def\N{\mathbb{N}}
\usepackage{ntheorem}

\newcommand{\calpha}{C$_\alpha$}

\newcommand{\VR}{\mathit{VR}}

\title{Homology of homologous knotted proteins}

\author[a]{Katherine Benjamin}
\author[a]{Lamisah Mukta} 
\author[a]{Gabriel Moryoussef}
\author[a]{Christopher Uren}
\author[a,b,1]{Heather A. Harrington}
\author[a,c,1]{Ulrike Tillmann}
\author[a,d,1]{Agnese Barbensi}

\affil[a]{Mathematical Institute, University of Oxford, Oxford OX2 6GG, United Kingdom}
\affil[b]{Wellcome Centre for Human Genetics,
University of Oxford OX3 7BN, United Kingdom}

\affil[c]{Isaac Newton Institute for Mathematical Sciences, University of Cambridge, CB3 0EH, United Kingdom}
\affil[d]{School of Mathematics and Statistics, University of Melbourne, Melbourne VIC 3010, Australia}

\affil[1]{To whom correspondence should be addressed:   E-mail: barbensi@maths.ox.ac.uk, tillmann@maths.ox.ac.uk, harrington@maths.ox.ac.uk}


\begin{abstract}
Quantification and classification of protein structures, such as knotted proteins, often requires noise-free and complete data. Here we develop a mathematical pipeline that systematically analyzes protein structures. We showcase this geometric framework on proteins forming open-ended trefoil knots, and we demonstrate that the mathematical tool, persistent homology, faithfully represents their structural homology. This topological pipeline identifies important geometric features of protein entanglement and clusters the space of trefoil proteins according to their depth.  Persistence landscapes quantify the topological difference between a family of knotted and unknotted proteins in the same structural homology class. This difference is localized and interpreted geometrically with recent advancements in systematic computation of homology generators. The topological and geometric quantification we find is robust to noisy input data, which demonstrates the potential of this approach in contexts where standard knot theoretic tools fail.

\end{abstract}

\begin{document}

\flushbottom
\maketitle

\thispagestyle{empty}

\section*{Introduction}
There are over 1,000 knotted protein structures currently cataloged,\cite{dabrowski2019knotprot} and the knotted domains in families of proteins with significant sequence differences have been conserved throughout evolution.\cite{sulkowska2012conservation}
The presence of knots in proteins substantially slows down the folding process,\cite{jackson2017fold, sulkowska2020folding} making the evolutionary selection somewhat counter-intuitive.
Mathematical and computational models predict that knots in proteins increase structural stability,\cite{mallam2010experimental,capraro2016untangling,dabrowski2016} which provides an explanation for this evolutionary selection. A clear experimental description of protein folding dynamics (\emph{e.g.,}~why and how knots arise) remains an open problem.\cite{jackson2020there,sulkowska2020folding}

Topological and geometric differences in knotted proteins correspond to differences in their folding pathways and dynamics as well as in their structural stability.\cite{dabrowski2016, dabrowski2019knotprot, barbensi2021depth,piejko2020folding,barbensi2021topological} In particular, the location and relative length of the entanglement (known as \emph{knot depth}) significantly influences the folding behavior of trefoil proteins.\cite{piejko2020folding} Perhaps most interesting is the case of knotted and unknotted carbamoyltransferases (AOTCases and OTCases), a pair of homologous proteins (\emph{i.e.} proteins with strong sequence and structural similarities), where a single crossing change which creates the entanglement is responsible for an increased stability of the knotted structures.\cite{sulkowska2008stabilizing}
To fully understand the interplay between entanglement, function, and folding, it is necessary to incorporate a nuanced topological and geometric characterization of protein structure alongside sequence and functional analysis. 

Recent studies have approached the characterization of knotted protein shape with techniques arising from low-dimensional topology and knot theory,\cite{tubiana2011probing, goundaroulis2020knotoids, barbensi2021depth} resulting in significant improvements in the classification of open knots and a better understanding of the knot folding mechanism.\cite{barbensi2021topological}
Although they are very powerful, low-dimensional topological techniques are computationally expensive and require precise and accurate input data.

Persistent homology (PH), the predominant tool in computational topology,\cite{zomorodian2005computing,ghrist2008barcodes,edelsbrunner2010computational} has enabled the characterization of meaningful topological and geometric features in data. PH quantifies features such as loops and voids in data at multiple scales of resolution, which can be taken as a fingerprint of the shape of data. Advances in PH, both in computational speed\cite{roadmap} and statistical tools,\cite{bubenik2015statistical,bubenik2017persistence} have enabled analysis of complex real-world data sets,\cite{vipond2021multiparameter,hiraoka2016hierarchical,mcguirl2020topological,stolz2020geometric} including protein structure and folding.\cite{xia2014persistent,gameiro2014topological,kovacev2016using,hamilton2021persistent}
A direct interpretation of homology features requires computing specific homology generators; their automatic and efficient computation is only recently possible.\cite{henselmanghristl6}

In this work we combine PH and low-dimensional topology to characterize geometric features of (open) knotted proteins. To our knowledge, the only other related work linking PH and low-dimensional topology is recent, and considers only closed knots.\cite{celoria2021statistical}
We propose and implement a computationally feasible PH pipeline for studying knotted proteins, at a global (full sequence) and local (substructural) scale, which does not rely on complex and computationally expensive knot invariants and sub-chain analysis.\cite{dabrowski2019knotprot, goundaroulis2020knotoids,barbensi2021depth,barbensi2021topological} Briefly, we input the 3D coordinates
of a knotted protein's structure, construct a point cloud spanning the backbone of the protein, compute the corresponding PH groups, and then translate these to persistence landscapes for protein comparison. We show that the information encoded in persistence diagrams and landscapes is enough to recover a clustering of the set of trefoil-knotted proteins by sequence similarity, thus validating PH as an effective fingerprint of protein structure. With this PH pipeline, we 
distinguish deeply versus shallowly knotted trefoil proteins and 
quantify structural differences between similar protein sequence homology classes. 
Computation of homology generators\cite{henselmanghristl6} has been used to show that the most persistent loops in PH correspond to active sites in the protein;\cite{kovacev2016using} here we show further that it distinguishes homologous knotted and unknotted proteins.

In particular, we are able to isolate the single local change that creates non-trivial entanglement in knotted AOTCases.
Topological statistics\cite{bubenik2015statistical,bubenik2017persistence}
allows us to demonstrate that these results are statistically significant. 
We showcase the robustness of this pipeline to incomplete or corrupt sequence data, which arises with many biopolymers.\cite{dabrowski2019knotprot, goundaroulis2020chromatin, siebert2017there} Hitherto, such analysis of noisy entangled curves is out of reach with a knot theoretic analysis.

\begin{figure}
\centering
\includegraphics[width=11cm]{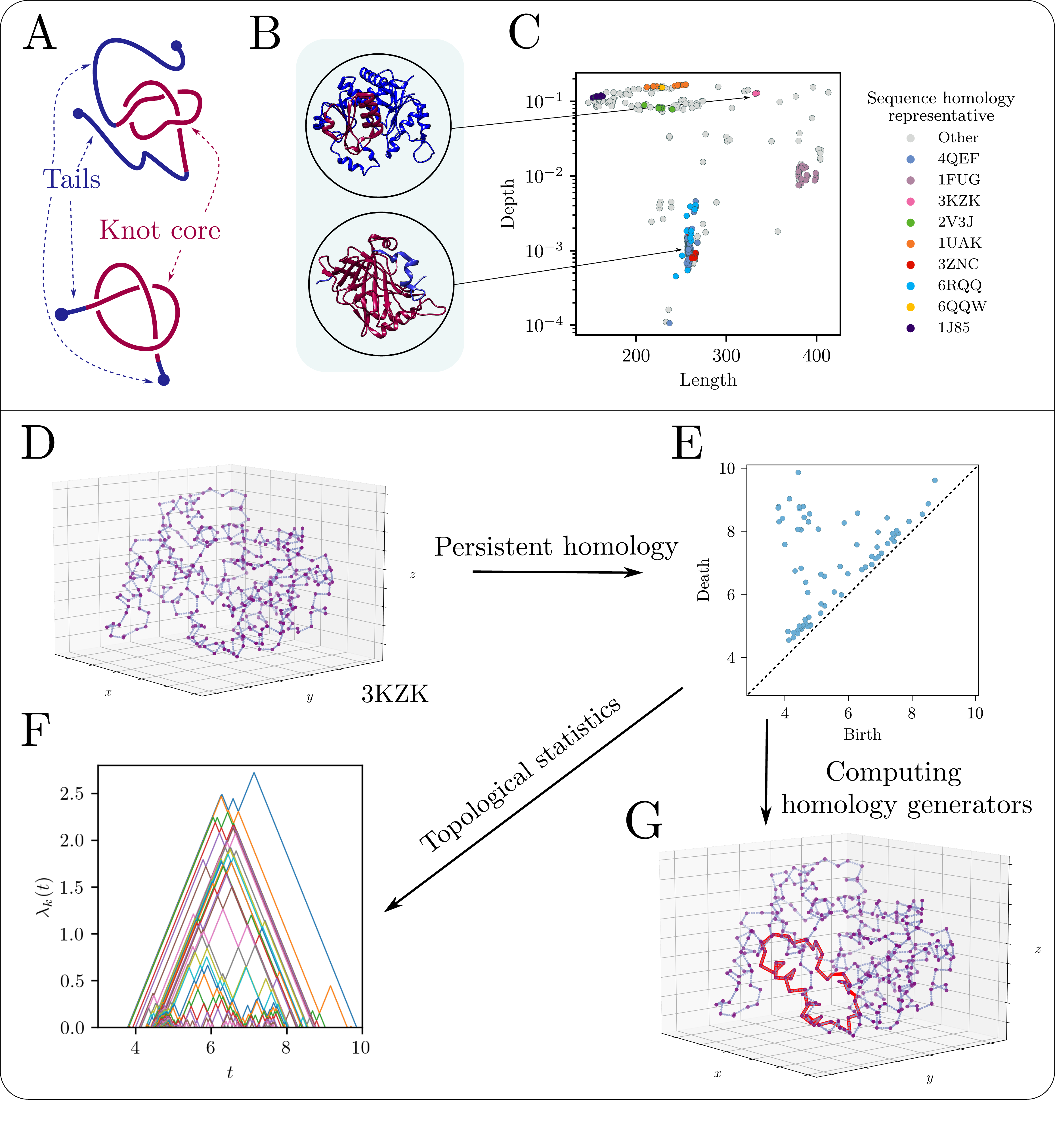}
\caption{\textbf{Data set and PH pipeline.} \textbf{(Data set)} \textbf{(A)} Schematics of a deeply knotted (top) and shallowly knotted (bottom) open curve. Knot cores tails are shown in purple and blue. \textbf{(B)} Example of a deeply knotted protein (PDB entry 3KZK) and a shallowly knotted one (PDB entry 4QEF). \textbf{(C)} The space of trefoil-knotted proteins plotted by chain length and knot depth. Each  protein is colored according to its sequence homology class.  Note that there are deeply and shallowly knotted proteins of the same length, as well as distinct sequence homology classes exhibiting similar length and depth. \textbf{(Pipeline)} \textbf{(D)} The protein data set is given by lists of 3D coordinates of \calpha{} atoms. For each protein, we generate the point cloud consisting of these points and linearly interpolated points between each successive \calpha{} atom.  \textbf{(E)} Persistence diagram derived from the 3KZK point cloud. The points represent 1-dimensional features corresponding to loops, and their positions represent the lifetimes of these features: their coordinates are their birth and death scales. \textbf{(F)} Persistence landscape derived from 3KZK. \textbf{(G)} PH generators in homology degree one can be represented by piece-wise linear cycles whose vertices are points in the point cloud. In red, an example of a local generator for a one-dimensional feature of the 3KZK point cloud.  }\label{fig:preliminaries}
\end{figure}

\section*{Overview of data set and pipeline}

\subsection*{Data set}
The first data set we consider consists of proteins whose backbones form an open-ended positive trefoil knot.\cite{dabrowski2019knotprot} The \emph{knot core} of a protein is defined as the shortest knotted sub-chain in its backbone. A protein is \emph{deeply knotted} if its knot core is entirely contained in a small portion of the chain, placed far away from the endpoints, see Figure \ref{fig:preliminaries}(A-C). Depth can be formally defined and quantified.\cite{barbensi2021depth} We label trefoil proteins as either \emph{shallow}, \emph{deep} or \emph{neither} based on their depth value. We also subdivide these trefoil-knotted proteins in the data set into structural homology classes, based on sequence similarity. Note that this division into classes can be performed using standard tools, see SI Appendix, Section 1.
The second data set consists of unknotted proteins sharing the same structural homology class as 3KZK, a deeply-knotted protein. This unknotted class is represented by 4JQO.

\subsection*{Pipeline}
For a given protein, we  first construct a point cloud  approximating its backbone curve, Figure \ref{fig:preliminaries}(D), and then  compute the PH of this point cloud in homology degree one. The result of this computation is summarized in the persistence diagram, Figure \ref{fig:preliminaries}(E), where each point with coordinates $(b,d)$ represents a 1-dimensional feature in the point cloud that is born at scale $b$ and dies at scale $d$. Next we translate the persistence diagram into the corresponding collection of landscapes; see Figure \ref{fig:preliminaries}(F). A point $(b,d)$ in the persistence diagram corresponds in the persistence landscape to a peak of height $(d-b)/2$ supported on the interval $[b,d]$. Persistence diagrams and persistence landscapes are two equivalent representations of PH. Next, for a given point $(b,d)$ of the persistence diagram or peak in the persistence landscape, we find a homology generator. A generator is represented by a sequence of sub-chains of the protein backbone which are linked end-to-end to form a loop, called a \emph{cycle}; see Figure \ref{fig:preliminaries}(G).

\subsection*{Experiments}

For the first experiment we consider two distinct measures of distance on PH: the Wasserstein distance  on persistence diagrams and the $L_1$ distance on persistence landscapes. 
Each distance induces a distinct dissimilarity measure on protein space, which we subsequently interpret using the dimensionality reduction algorithm Isomap.\cite{tenenbaum2000}
We infer the typical PH of each structural homology class by computing its average persistence landscape.

For the second experiment, we consider a collection of trefoil-knotted and unknotted proteins belonging to the same structural homology class. We compute a generator of a 1-dimensional feature which is unique to the PH of the knotted class.

In the third experiment, we reproduce the results of the first two experiments (a) with a sparser sampling of the protein backbone and (b) after perturbing each point cloud by increasing amounts of Gaussian noise.

\section*{Results}
\subsection*{Global clustering and depth type detection}

We demonstrate that PH captures structural, geometric and topological differences in proteins whose backbones form trefoil knots.

\begin{figure}
\centering
 \includegraphics[width=10cm]{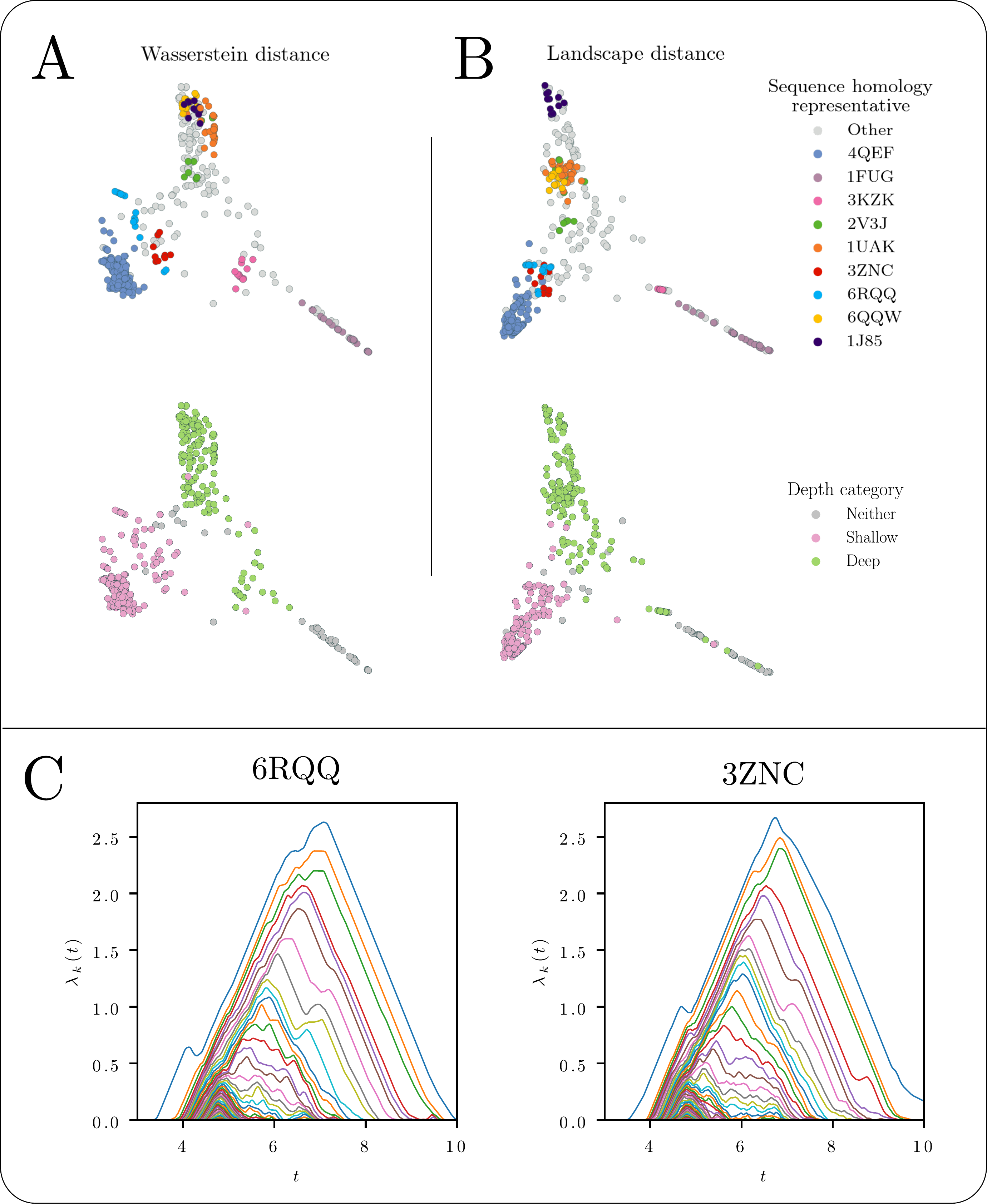} \caption{\textbf{Global analysis: the space of knotted protein structures.}  \textbf{(A)} Isomap embedding of the space of trefoil-knotted proteins equipped with the Wasserstein distance on persistence diagrams (see SI Appendix, Section 2). Given a distance matrix, Isomap produces a configuration, 2-dimensional in our case, such that the new distance between any two objects is preserved as much as possible. The embedding forms clusters corresponding to sequence homology classes. The embedding successfully clusters by depth category. \textbf{(B)} Isomap embedding of the space of trefoil-knotted proteins equipped with the distance on persistence landscape. For a definition of this distance see SI Appendix, Section 2. The embedding forms clusters corresponding to sequence homology classes. The embedding successfully clusters by depth category. \textbf{(C)} Average persistence landscapes generated from the sequence homology classes with representative 6RQQ (left) and 3ZNC (right). Although these classes are not separated in the Isomap embeddings in (A) and (B), a randomization test confirms that the difference in their average landscapes is statistically significant ($p\approx0.003$). }\label{fig:2}
\end{figure}

First we compute the persistence diagrams and landscapes in dimension one for every protein in the data set, and then pairwise distances between persistence diagrams and landscapes. 
The resulting structures induced by these distances on the space of trefoil proteins are approximated by dimensionality reduction, shown in Figure \ref{fig:2}(A-B, top panel).

Proteins that are in the same homology class cluster together by topological distance. This clustering validates that the PH pipeline recovers structural features of knotted proteins. We also observe separation between distinct homology classes. In particular, those sequence homology classes exhibiting significantly different chain lengths or depth types are very well separated, for example 1FUG and 3KZK. 
While some of the largest structural homology classes overlap in the projections, taking both the Wasserstein and $L_1$ landscape metrics fully distinguishes the different sequence homology classes. We emphasize that the two distance matrices and corresponding Isomap projections are computed from the same initial data: the PH of each protein point cloud. 

We apply topological statistics to infer the PH shape of structural homology classes, and detect geometric differences between proteins in these classes. 
We define the topological fingerprint of each structural homology class to be the average landscape of the class (see Figure \ref{fig:2}(C)).
We employed a 1000-sample randomization test to compare these topological fingerprints
(see SI Appendix, Section 2 and Table S1). The test produced an approximate $p$-value of 0.001 for all pairs of classes except 3ZNC and 6RQQ, which gave a $p$-value of $0.003$. 
Therefore, PH detects statistically significant pairwise differences between all structural homology classes, including those not separated in the Isomap configurations.

Using this topological pipeline we next detect knot depth, which is a meaningful geometric feature intrinsic to protein entanglement. Figures \ref{fig:2}(A-B, bottom panel) show a considerable separation between shallowly and deeply knotted proteins. In particular, even deeply knotted proteins and shallowly knotted proteins having similar chain length, as those in the 4QEF and 1UAK classes, are clearly separated. This shows that, beyond separating proteins by structural similarity, PH can be successfully employed to analyze the geometric entanglement of open curves. 

\subsection*{Detecting local topological changes in homologous proteins}
\begin{figure}
\centering
\includegraphics[width=10cm]{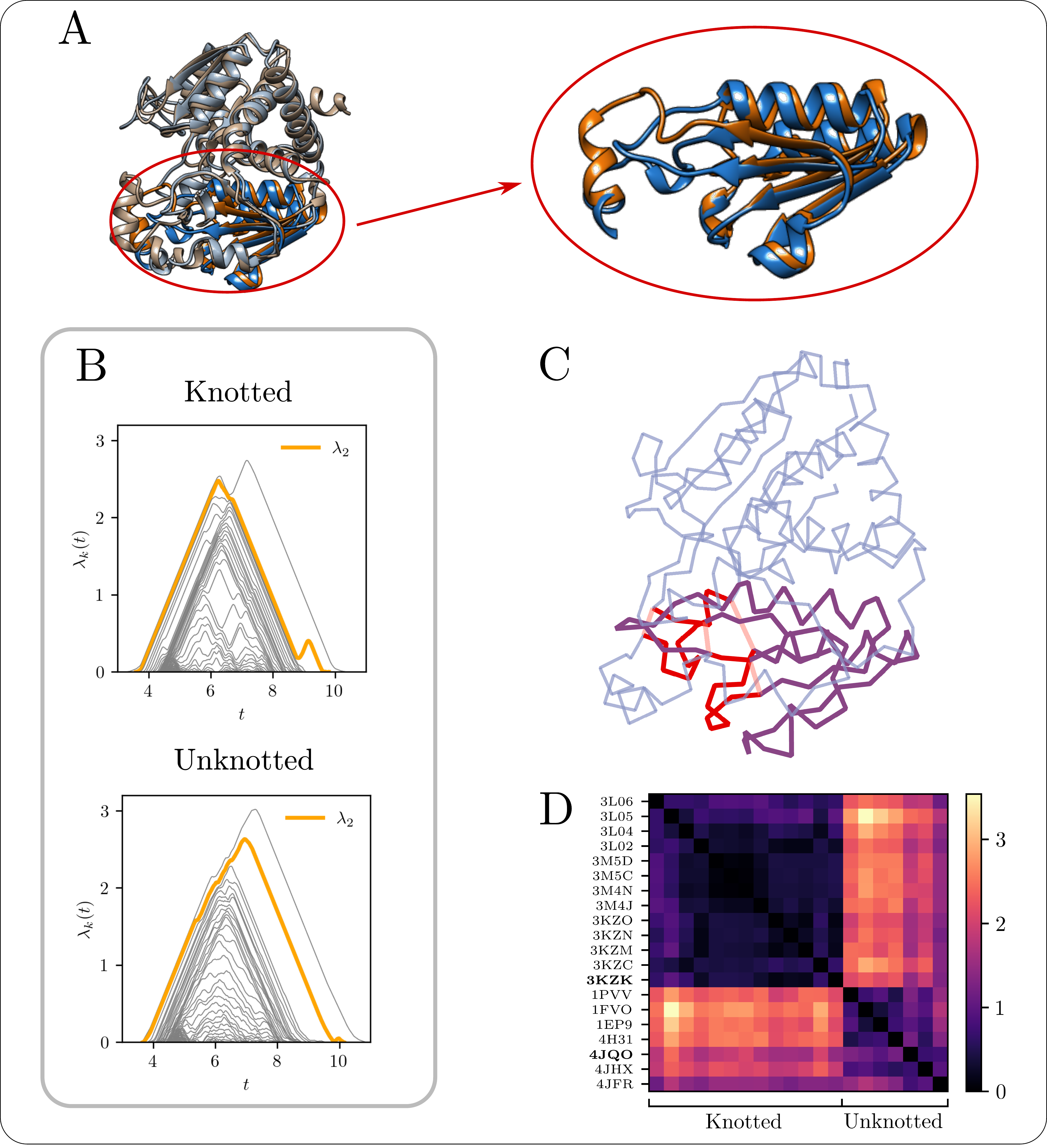}
\caption{\textbf{Local analysis: geometry of homologous protein substructures.} \textbf{(A)} Two homologous proteins, 3KZK (blue, knotted) and 4JQO (orange, unknotted), overlaid. These proteins have almost superimposable structures, but differ as knots by a crossing change localized within the red ellipse. The knot core in 3KZK and its corresponding structure in 4JQO are highlighted by showing the remaining parts in lighter shades of orange and blue. On the right-hand side, a close up of the local configurations causing the topological change. A strand movement transforms the deeply knotted 3KZK into the unknotted 4JQO. \textbf{(B)} Average persistence landscapes generated from knotted (top) and unknotted (bottom) protein chains. The peak in $\lambda_2$ (orange) centered at $t \approx 9$ in the knotted case corresponds to a generator $c$ for the persistent homology of the knotted chains which does not arise in the persistent homology of the unknotted chains. \textbf{(C)} The backbone of 3KZK. Violet segments indicate the knot core. The cycle representing the persistent homology generator $c$ is plotted in red and pink, where the pink segments show simplices in $c$ that are not part of the 3KZK curve. Note that $c$ is positioned close to the knot core, and more specifically, close to the crossings responsible for the non-trivial entanglement. Further, $c$ intersects the arc that needs to be pushed to untangle the curve.  \textbf{(D)} Heat-map showing the distances between the $\lambda_2$ landscapes for the proteins in the AOTCase and OTCase families. The two distinct purple squares demonstrate sufficient similarity in each class for the average $\lambda_2$ landscapes to be faithful representatives for each class.}\label{fig:3} 
\end{figure}

We now investigate whether PH is sensitive to knottiness.
We consider the homologous pair of protein structures: the trefoil-knotted 3KZK and the unknotted 4JQO. As shown in Figure \ref{fig:3}(A), the difference in structure between these proteins is localized in a few separate non-overlaying portions. One of these portions encloses the crossing change responsible for the change in the knot type, and a natural question is whether PH is capable of detecting and localizing this topological difference. To this end, we extend the data set to all the trefoil-knotted AOTCases (of which 3KZK can be taken as a representative), and the unknotted OTCases (of which 4JQO can be taken as a representative).\cite{potestio2010knotted, dabrowski2016}

The average persistence landscapes for these two classes are displayed in Figure \ref{fig:3}(B). Despite the high similarity of all the protein structures considered here, we notice remarkable differences between the two average landscapes, most notably in the $\lambda_2$ landscape layers (Figure \ref{fig:3}(B), orange curve).

For the trefoil-knotted AOTCase family, $\lambda_2$ has a small second peak centered at scale value 9, which is not present for the unknotted OTCase family. To investigate this difference, we compare only the $\lambda_2$ layer for each of the proteins in the two classes. The heat-map shown in Figure \ref{fig:3}(D) confirms that the $\lambda_2$ layers within each of the AOTCase and OTCase families are sufficiently similar to each other, thus the average $\lambda_2$ landscape layer is a good representative of each class.

We geometrically interpret the $\lambda_2$ layer difference by analyzing the corresponding homology generators. We focus on the $\lambda_2$ layers of specific proteins, namely the knotted AOTCase 3KZK and the OTCase 4JQO. In the landscape of 3KZK, the $\lambda_2$ peak corresponds to a specific PH generator $c$ via the landscapes-diagrams correspondence.\cite{bubenik2015statistical} A cycle representing the PH generator $c$, computed using a recently available and computationally efficient algorithm in Eirene\cite{henselmanghristl6} (see Figure \ref{fig:3}(C), red). Most of the simplices forming the cycle are part of the  knot core (Figure \ref{fig:3}(C), violet) of the protein backbone. The homology generator identifies the portion of the protein backbone that forms the protein knot. Specifically, the cycle overlaps with the essential crossing that distinguishes the knotted and unknotted proteins.
We compute and find a similar result across all of  the AOTC structures. 

\subsection*{The PH pipeline is robust}
One of the main strengths of PH is its theoretically-guaranteed robustness to noise.\cite{cohen2007stability} A natural question is whether the results stand when the input data is incomplete or noisy. To check this, we  first (a) apply the pipeline to a sparse point cloud, given directly by the \calpha{} atoms for each protein -- that is, without interpolating between adjacent atoms.  We then (b) perturb each point cloud by applying increasing amounts of Gaussian noise to the coordinates of each point, and repeat the procedure yet again (see SI Appendix, Section 4). Remarkably, we can successfully reproduce both the global and local results of the previous sections. The sequence homology clustering only fails when perturbing the points with a standard deviation of approximately one third the distance between two consecutive amino-acids, while the crossing change detection is lost slightly earlier.
Isomap plots and the corresponding local generator analysis are described in detail in SI Appendix, Section 3 and Figures S1 and S2.

\section*{Discussion}
Our investigation illustrates the potential of using PH to analyze biopolymers with complex geometry. We showed Wasserstein distance and landscape distance meaningful discriminated between protein structures without the computationally expensive pre-processing, such as sequence alignment, required for traditional methods.

While structural classification of proteins via PH has previously been explored in other works,\cite{xia2014persistent,gameiro2014topological,kovacev2016using,hamilton2021persistent} only two have used the backbone curve as input data,\cite{xia2014persistent, hamilton2021persistent} and only the developer of landscapes used the statistical apparatus of persistence landscapes for this task.\cite{kovacev2016using} We successfully employed persistence landscapes to distinguish between structurally similar knotted and unknotted proteins, thus detecting their global topology. Furthermore, we used generators to localize the differences between the topological types of homologous proteins. This proposed methodology has the potential to be applicable in other contexts in which knots arise naturally.

Here we built, combined and integrated on all of these techniques to explore the space of knotted proteins and to detect features specific to open knots, for example knot depth. Our results demonstrate that PH is sensitive to topological and geometric features of those structures. Furthermore, we demonstrated that our approach works with noisy data and could therefore be a computationally efficient tool for the study of knotted structures in other contexts. 
The robustness of our methodology could be particularly crucial in cases in which experiments are performed at a resolution that does not guarantee a complete determination of a biopolymer's underlying spatial curve, as for example in the case of DNA chromosomes.\cite{goundaroulis2020chromatin, siebert2017there}

\section*{Materials and Methods}

\subsection*{Protein data}
We model each protein structure as the piece-wise linear (PL) curve spanned by its \calpha{} atoms, and we consider proteins whose backbones form an open-ended positive trefoil.\cite{dabrowski2019knotprot} We construct a point cloud from its PL backbone curve by linearly interpolating a total of 5 equidistant points between each pair of consecutive \calpha{} atoms in the curve. For more details, see SI Appendix, Section 1.

\subsection*{PH computations}
Persistence diagram in dimension one are computed using Ripser.\cite{bauer2021ripser} Pairwise $W_1[L_\infty]$-Wasserstein distances are computed using GUDHI.\cite{gudhi} Persistent homology generators are computed using Eirene.\cite{henselmanghristl6} Persistence landscapes are computed using Python software, available at \url{https://github.com/katherine-benjamin/ph-knotted-proteins}. Our software is based on the algorithms given in.\cite{bubenik2017persistence} The software includes scripts for finding and plotting landscapes from a barcode, computing average landscapes, and calculating distance. For more details, see SI Appendix, Sections 2 and 4.

\subsection*{Noisy data}
To probe our methodology's robustness we apply our pipeline to the non-interpolated point clouds obtained by considering only the \calpha{} atoms, and to successive perturbations of such point clouds.
For more details, see SI Appendix, Sections 2 and 4 and Figures S1 and S2. 

\subsection*{Data Availability}
The data and code are available on the GitHub repository: \\ \url{https://github.com/katherine-benjamin/ph-knotted-proteins}.

\section*{Acknowledgements}
The authors thank Greg Henselman for suggesting to analyze generators. AB, KB, HAH, and UT are grateful to the support provided by the UK Centre for Topological Data Analysis EPSRC grant EP/R018472/1. HAH gratefully acknowledges funding from EPSRC EP/R005125/1 and EP/T001968/1, the Royal Society RGF$\backslash$EA$\backslash$201074 and UF150238.

\section*{Author contributions statement}
Author contributions: A.B., H.A.H., and U.T. conceived the research; K.B., A.B., H.A.H., and U.T. designed the research; K.B., L.M., G.M., and C.U. performed research; K.B., L.M., G.M., C.U., and A.B. analyzed data; and all the authors wrote the paper.

\section*{Additional information}
The authors declare no competing interest.

\newpage
\begin{appendices}
These appendices contain the Supplementary Information for: ``Homology of homologous knotted proteins''.

\section{Data}

\subsection{Protein structure}

Proteins are long chains of amino acids that fold into specific three-dimensional structures.
Proteins fold into complex 3D structures according to this sequence, and proteins with sufficiently similar sequences of amino acids are highly likely to have similar 3D structures. Such proteins are called \emph{homologous}.

Each amino acid in a protein chain contains a distinguished atom, the \calpha{} atom, which can be taken as its representative in space. 
For solved protein structures, the spatial coordinates of \calpha{} atoms are recorded and stored in the Protein Data Bank (PDB).\cite{berman2007worldwide}
Protein structures are often modeled\cite{dabrowski2019knotprot,goundaroulis2020knotoids,barbensi2021depth,barbensi2021topological} as the piece-wise linear (PL) curve,
called the \emph{protein backbone},
spanned by their \calpha{} atoms. In this model, the \emph{length} of a protein is given as the number of amino acids in its sequence (\emph{i.e.}~as the number of vertices in the corresponding PL curve).

\subsection{Knotted proteins}

Some protein backbones are entangled in a non-trivial way.\cite{dabrowski2019knotprot} Detecting and characterising entanglement for abstract open curves is a non-trivial task, since, differently from the case of closed curves,\cite{rolfsen2003knots} they do not have an \emph{a priori} well defined knot type. In order to identify the topological type of an open curve, one can construct a closed curve from the initial open one by joining the two endpoints with a new segment. 
The resulting closed curve has a knot type that can be detected using knot invariants. 
This knot type, however, is dependent on the additional segment chosen. 
To overcome this issue, the most common approaches use probabilistic closures of the open chain, and then take the most likely knot type obtained to be the knot type of the initial open curve.\cite{dorier2021open} When this knot type is non-trivial, we say that the curve forms an open-ended knot.

Similarly, it is possible to detect the \emph{knot core} of a curve, intuitively defined as its shortest knotted sub-chain. This can be done by performing the so-called \emph{sub-chain analysis}, in which the knot types of all the sub-chains of the curve are computed.\cite{king2007identification,millett2013identifying, dabrowski2019knotprot,dorier2018knoto} There are several different ways to compute the location of a curve's knot core; in this paper, we use the results stored in the online database KnotProt,\cite{dabrowski2019knotprot} where knot types and knot cores are computed using the software Knoto-ID.\cite{dorier2018knoto}

We can therefore partition a knotted open curve into the knot core and two tails: the N-tail and the C-tail. For a given protein chain $p$, the \emph{knot depth}\cite{barbensi2021depth} is defined as
\begin{equation}
    D(p) = \frac{l(N)l(C)}{l(T)^2},
\end{equation}
where $N$, $C$, and $T$ refer to the N-tail, C-tail, and the whole curve respectively, and the length function $l$ counts the number of \calpha{} atoms in a subchain.

We distinguish between three types of knotted proteins depending on their depth: Proteins with depth above $0.05$ are referred to as
\emph{deeply knotted}; those with depth below $0.005$  as \emph{shallowly knotted}; and  other proteins as \emph{neither}. 

\subsection{Data set}

The data set for this paper is drawn from two sources. The first is a list of 513 trefoil-knotted protein chains, which comprises nearly all\footnote{We excluded 4 chains with particularly high backbone length (over 500 \calpha{} atoms), and the remaining 20 missing chains were cataloged after we performed our analysis.} of the 537 trefoil-knotted chains catalogued to date in KnotProt.\cite{dabrowski2019knotprot} For each such chain, we also take its length and the lengths of its tails from KnotProt, and use this to compute the knot depth.
We only consider proteins whose chains form open ended trefoils, as these are by far the largest group of knotted proteins (the database KnotProt currently contains 668 different knotted chains). Moreover, for proteins having different knot types, there is not substantial variation of depth and structural homology type.

We divide the set of trefoil-knotted chains into sequence homology classes, according to the procedure in.\cite{barbensi2021topological} In detail, we first perform a sequence similarity search using the RCSB PBD Search API,\cite{rose2021} and construct a sequence similarity matrix from the results. From this, we perform single-linkage hierarchical clustering using SciPy,\cite{2020SciPy-NMeth} and then impose a 70\% distance threshold to obtain the different clusters. We keep the largest 9 clusters as structural homology classes and label those by the PBD ID of a representative, and label everything else as `Other'. Information on knotting and sequence homology for each trefoil-knotted chain that we use is stored in the GitHub repository.\cite{ph-knotted-proteins}

For the second source of data set, we consider the AOTCase and OTCase protein structures, which are respectively knotted and unknotted. The list of these structures, consisting of 13 knotted chains and 7 unknotted chains, can be found in the GitHub repository.\cite{ph-knotted-proteins}

For every protein chain considered, we store the 3D \calpha{} coordinates in \texttt{.xyz} files. Individual lines contain the 3D coordinates of successive \calpha{} atoms in the protein backbone. For knotted proteins, we download these directly from KnotProt.\cite{dabrowski2019knotprot} For unknotted proteins, we first download protein data from the Protein Data Bank (PDB),\cite{berman2007worldwide} and then convert them to this standard format using the \texttt{pdb\_to\_xyz.R} script provided as part of Knoto-ID 1.3.0.\cite{dorier2018knoto}

\subsection{Preprocessing} 
In order to compute PH, we need a point cloud which approximates the PL curve representing the protein chain. To obtain such a point cloud, we preprocess protein data by adding interpolated points between successive \calpha{} atoms along the backbone of each structure.  In particular, suppose $p_1, p_2 \in \mathbb{R}^3$ are the coordinates of successive \calpha{} atoms in a protein chain. We add $d$ additional equidistant, linearly interpolated points $x_1, \dots, x_d$ according to the formula
\begin{equation}
    x_i = p_1 + i\left(p_2 - p_1\right)/(d+1).
\end{equation}

We then repeat this for each successive pair of \calpha{} atoms to arrive at an interpolated protein chain, which we treat as a point cloud and use as input for the persistent homology calculations. We chose $d=5$ for our analysis in order to strike a balance between geometric faithfulness and computational feasibility.

\section{Robustness results}
We repeat all the computations and reproduce all of our results by taking as input data the point cloud given by the \calpha{} atoms (\emph{i.e.}~, without interpolating). The results are shown and discussed in Figure \ref{fig1}.

We further repeat all the computations after adding an increasing amount of Gaussian noise to each point cloud (see Section \ref{sec:noise}). The results are shown and discussed in Figure \ref{fig2}.

\section{Topological methods}

\subsection{Persistent homology} 

\paragraph{Vietoris-Rips complex} Topological data analysis (TDA) is a field which aims to analyse the shape of data sets. To do this we want to build a suitable topological space from the data, and there are many ways to do this. In this application, for each protein we have constructed a point cloud in $\mathbb R^3$, which may be thought of as a metric space under the usual Euclidean distance, and the topological space we build from it is a simplicial complex, known as the \emph{Vietoris-Rips} complex.
Given a finite metric space $(M,d)$ and a non-negative real number $\epsilon$, the \emph{Vietoris-Rips} complex $\VR_\epsilon(M)$ at scale $\epsilon$ is
\begin{equation}
    \VR_\epsilon(M) = \{\sigma \subseteq M \mid d(x,y)\leq \epsilon\textrm{ for all } x,y \in \sigma\}.
\end{equation}

We recall the definition of a \emph{filtration}. An $\R^{\geq 0}$-indexed filtration is a topological space $X$ and a collection $\{X_t\}$ of subspaces of $X$ such that $X_t \subseteq X_s$ whenever $t \leq s$, and $X = \cup_{t \in \R^{\geq 0}} X_t$. One can easily check that by considering all possible values for the scale parameter $\epsilon$, the associated Vietoris-Rips complexes form an $\R^{\geq 0}$-indexed filtration.

\paragraph{Persistence} Persistent homology, a key tool from TDA, seeks to understand a data set by considering topological features that exist at a range of scales. One way to do this is to construct the filtered Vietoris-Rips complex of the data set and then consider the homology of each complex in the filtration.  For each homological degree $i\geq 0$, this yields a collection of vector spaces $H_i(\VR_\epsilon(M))$ indexed by the scale parameter $\epsilon$ and linked by linear maps $H_i(\VR_\epsilon(M)) \to H_i(\VR_\delta(M))$ induced by the inclusion $\VR_\epsilon(M) \subseteq \VR_\delta(M)$ for each $\epsilon \leq \delta$. As the scale parameter $\epsilon$ increases, different topological features appear in the complexes $\VR_\epsilon(M)$, and these appear as generators for the vector spaces $H_i(\VR_\epsilon(M))$. We can record the lifespans of these features in a \emph{persistence diagram}: if a feature appears in the homology at scale value $b$ and disappears at scale value $d>b$, we record this with the pair $(b,d)$. The persistence diagram of the filtration is the multi-set of all such pairs, known as birth-death pairs.

The robustness of PH to noise is guaranteed by the Stability Theorem\cite{cohen2007stability} which states that similar input data will produce similar persistence diagrams. This is crucially important for biological applications as noise is inevitably introduced during experimental data collection, as well as the inherent noise that exists in, for example, in vivo protein configurations.

\paragraph{Generators}
We are most interested in the first homology, $H_1(VR_\epsilon (M))$, and will discuss particular elements. A generator (or basis element) of this vector space is given by a 1-dimensional {\it cycle}, that is a sequence  $(m_1, \dots, m_k) $  of points in $M$ 
with $d(m_i, m_{i+1}) \leq \epsilon $ and $d(m_1, m_k) \leq \epsilon$. Any such cycle represents a generator if it is not a {\it boundary} of a 2-dimensional chain. Intuitively this means that the 1-dimensional 'hole' defined by the cycle is not filled in by 2-simplices. There are in general many cycles representing the same generator in homology. Any two such choices differ by a boundary. For more details, see.\cite{hatcher}

\paragraph{Distances}
To compare the results of persistent homology, we define a distance metric on the space of persistence diagrams, and the most commonly used are the $L_q$-Wasserstein distances, defined for diagrams $D^1$ and $D^2$, as
\begin{equation}
W_p[L_q](D^1,D^2) = \operatorname{inf}_{\phi:D^1\to D^2}\left [ \sum_{I\in D^1} L_q[I,\phi(I)]^p \right ] ^{1/p}
\end{equation}
 where $\phi$ ranges over all partial maps from $D^1$ to $D^2$.

In our analysis we take $p=1$ and $q=\infty$, and refer to the $W_1[L_\infty]$-Wasserstein distance simply as Wasserstein distance. 

\subsection{Persistence landscapes} \label{sec:landscapes}

\paragraph{Construction} Persistence landscapes are a vectorization of persistence diagrams which enable the use of a wide variety of statistical tools.\cite{bubenik2015statistical} To construct them, we consider a point $I = (b,d)$ in a persistence diagram $D$. We define the auxiliary function
\begin{equation}
    f_I(t) = \begin{cases} 
0 & \text{if } t \notin [b,d], \\
t-b & \text{if } t \in [b, \frac{b+d}{2}], \\
-t+d & \text{if } t \in [\frac{b+d}{2},d].
\end{cases}
\end{equation}
We then define the \emph{persistence landscape corresponding to $D$} as the function $\lambda \colon \N \times \R \to \R$ where $\lambda(k,t)$ is the $k$-th largest value of $f_I(t)$ over each $I \in D$. Equivalently, we view a persistence landscape as being given by a sequence of functions $\lambda_k \colon \R \to \R$ given by $\lambda_k(t) = \lambda(k,t)$.

\paragraph{Distances} There are several notions of distance or dissimilarity that can be defined on the space of persistence landscapes. The standard distance is derived from the $p$-norm $\lVert \lambda \rVert_p$ given by
\begin{equation}
    \lVert \lambda \rVert_p^p = \sum_{k \geq 0} \int_\R \lvert \lambda_k \rvert^p.
\end{equation}
In particular, we then derive an $L_p$ distance between two landscapes $\lambda^1$ and $\lambda^2$ given by the $p$-norm of their pointwise difference:
\begin{equation}
    d_p(\lambda^1, \lambda^2) = \lVert \lambda^1 - \lambda^2 \rVert_p.
\end{equation}
We refer to $d_1$ as the \emph{$L_1$ distance on persistence landscapes}.

In addition, it is often useful to compare individual layers of persistence landscapes. In particular, if $S \subseteq \N$, we have a difference measure
\begin{equation}
    d^S(\lambda^1 - \lambda^2) = \sum_{k \in S} \int_\R \lvert \lambda_k^1 - \lambda_k^2 \rvert
\end{equation}
comparing just the layers $\lambda_k$ for $k \in S$. Bubenik proves that, under some easily-satisfied boundedness conditions, this structure enjoys several useful statistical properties.\cite{bubenik2015statistical} Note that when comparing knotted AOTCases and unknotted OTCases we consider only the case when $S = \{2\}$, since the main difference in the average landscapes was contained there, as explained in the Results section of the main manuscript.

\paragraph{Averages} Given a collection $\lambda^1, \dots, \lambda^N$ of persistence landscapes, their average persistence landscape $\overline{\lambda}$ is defined pointwise as
\begin{equation}
    \overline{\lambda}(k, t) = \frac{1}{N} \sum_{i=1}^N \lambda^i(k,t).
\end{equation}
This is the landscapes' unique Fr\'echet mean.
\paragraph{Statistical tests}
To compare samples of persistence landscapes, we employ randomization tests.\cite{bubenik2017persistence} If $\lambda^1, \dots, \lambda^N$ and $\lambda'^1, \dots, \lambda'^M$ are two samples, we compare them with the test statistic
\begin{equation}\label{eq:teststatistic}
    t = \lVert\overline{\lambda} - \overline{\lambda'}\rVert.
\end{equation}
We can consider the samples as giving a partition of the set of $N + M$ landscapes into two groups of size $N$ and $M$.
To perform a $k$-sample randomization test, we produce $k$ distinct random partitions of the two samples into groups of size $N$ and $M$. 
The $p$-value produced is then the proportion of partitions which produce a test statistic at least as large as the test statistic produced from the original partition.

\section{Computation}

Python code for each stage of the computational pipeline is provided in the accompanying GitHub repository,\cite{ph-knotted-proteins} along with more detailed instructions on its use. In particular, it is possible to use this code to reproduce any programmatically generated plots displayed in the paper.

\subsection{Persistence diagrams}

We use two different software packages to compute the first persistent homology of the Vietoris-Rips complexes associated to interpolated protein chain structures. For most of the structures we do not compute representative generators. In these cases, we use Ripser 1.2.1\cite{bauer2021ripser} to generate persistence diagrams. Otherwise, for AOTCase and OTCase structures, we use Eirene 1.3.6\cite{henselmanghristl6} to produce both persistence diagrams and canonical representative generators for homology.
Once persistence diagrams have been computed, we use GUDHI 3.4.1\cite{gudhi} to compute pairwise Wasserstein distances between them.

\subsection{Persistence landscapes} \label{sec:landscapecomps}

To compute and manipulate persistence landscapes, we use a modified version of the Pysistence Landscapes software package (\url{https://gitlab.com/kfbenjamin/pysistence-landscapes/}) maintained by the first author. 
 This package is based on the algorithms described by Bubenik and Pawe\l{}.\cite{bubenik2017persistence} The modified version of the package is included in the GitHub repository.\cite{ph-knotted-proteins} Persistence landscapes are stored in human-readable \texttt{.lan} files as specified in.\cite{bubenik2017persistence} For statistical tests, we use the Pysistence Landscapes software to carry out 1000-sample randomisation tests according to the test statistic in \eqref{eq:teststatistic}. The results of these tests are recorded in Table \ref{tab:pvalues}.

\subsection{Dimensionality reduction}

We use scikit-learn 0.24.2,\cite{scikit-learn} with default settings, to compute two-dimensional Isomap embeddings based on the precomputed Wasserstein and persistence-landscape distance matrices.

\subsection{Gaussian noise}\label{sec:noise}

We simulate low resolution data by perturbing each point cloud with Gaussian noise. In particular, given a standard deviation $\sigma$ we replace each point $(x_1,x_2,x_3)$ in each point cloud with a `noisy' point
\[ 
 (x_1 + \varepsilon_1,  x_2 + \varepsilon_2, x_3+\varepsilon_3)
\]
where each $\varepsilon_i$ is drawn independently from $\mathcal{N}(0, \sigma^2)$. The Python code for this process is provided in the GitHub repository.\cite{ph-knotted-proteins}

We ran the experiments with standard deviations $\sigma \in \{0.1, 0.2, \dots, 1.0\}$ and the resulting outputs and figures are available in the GitHub repository.\cite{ph-knotted-proteins}

\subsection{Silhouette coefficient} \label{sec:silhouette}

We use silhouette coefficients\cite{rousseeuw1987silhouettes} to measure the faithfulness of Isomap embeddings as noise increases. These are computed using scikit-learn 0.24.2.\cite{scikit-learn}

\iftrue
\begin{figure}
\centering
\includegraphics[width=15cm]{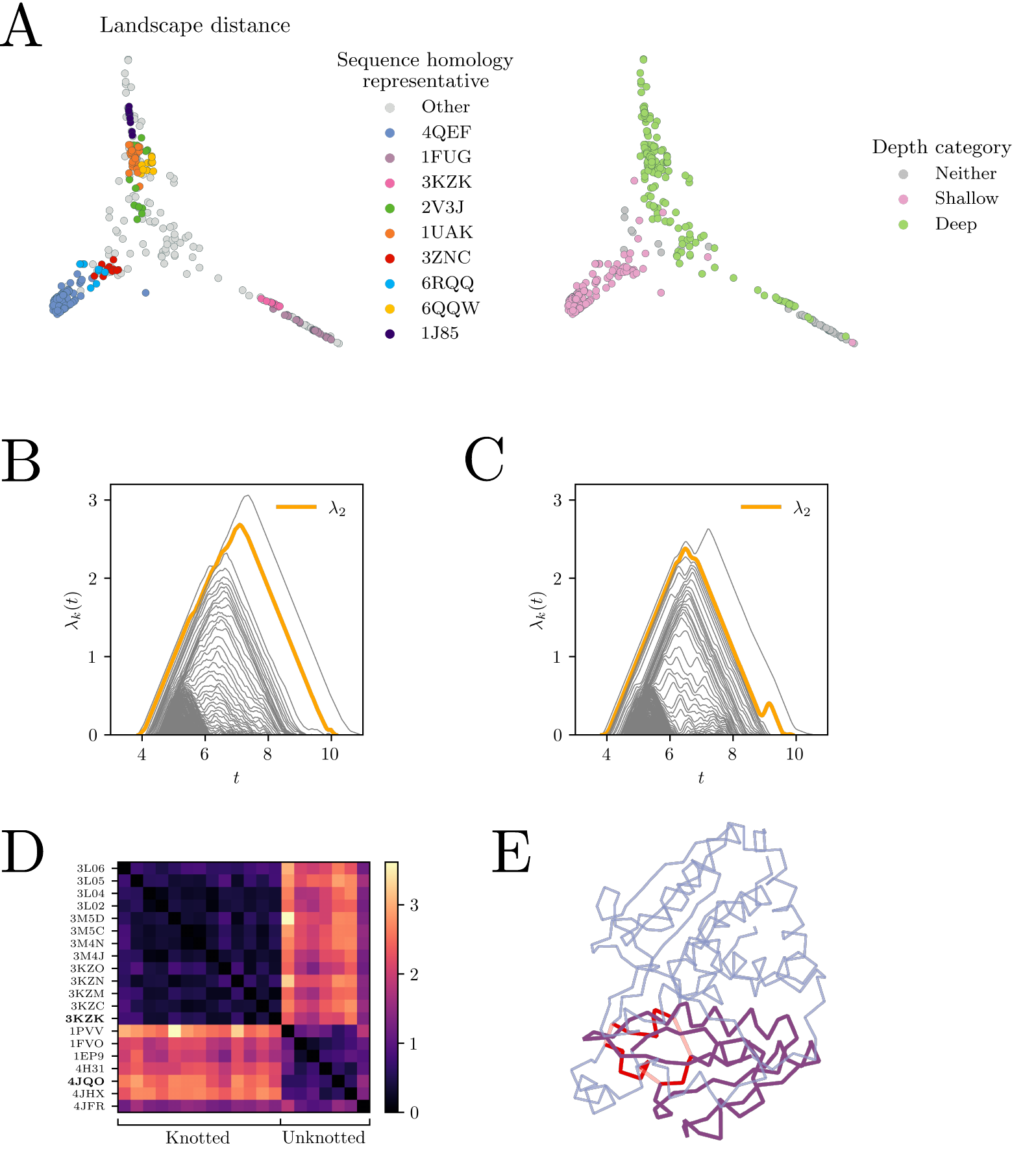}
\caption{\textbf{Sparse dataset} \textbf{(A)} Isomap embedding of the space of trefoil-knotted proteins equipped with the distance on persistence landscape.
The embedding still forms clusters corresponding to sequence homology and the depth types are still correctly detected. \textbf{(B)} The average landscape for the unknotted OTCases. \textbf{(C)} The average landscape for the knotted AOTCases. \textbf{(D)}  Heat-map showing the distances between the $\lambda_2$ landscapes for the proteins in the AOTCase and OTCase families. The two distinct purple squares still demonstrate sufficient similarity in each class for the average $\lambda_2$ landscapes to be faithful representatives for each class. \textbf{(E)} The backbone of 3KZK. Violet segments indicate the knot core. The peak in $\lambda_2$ (orange) centered at $t \approx 9$ in the knotted case (C) corresponds to a generator $c$ for the persistent homology of the knotted chains which does not arise in the persistent homology of the unknotted chains (B). The cycle representing $c$ is plotted in red and pink, where pink segments are simplicises in $c$ that are not part of the 3KZK PL curve. Note that $c$ is still positioned close to the knot core, and more specifically, close to the crossings responsible for the non-trivial entanglement. Further, $c$ still intersects the arc that needs to be pushed to untangle the curve.}
\label{fig1}
\end{figure}

\begin{figure}
\centering
\includegraphics[width=17cm]{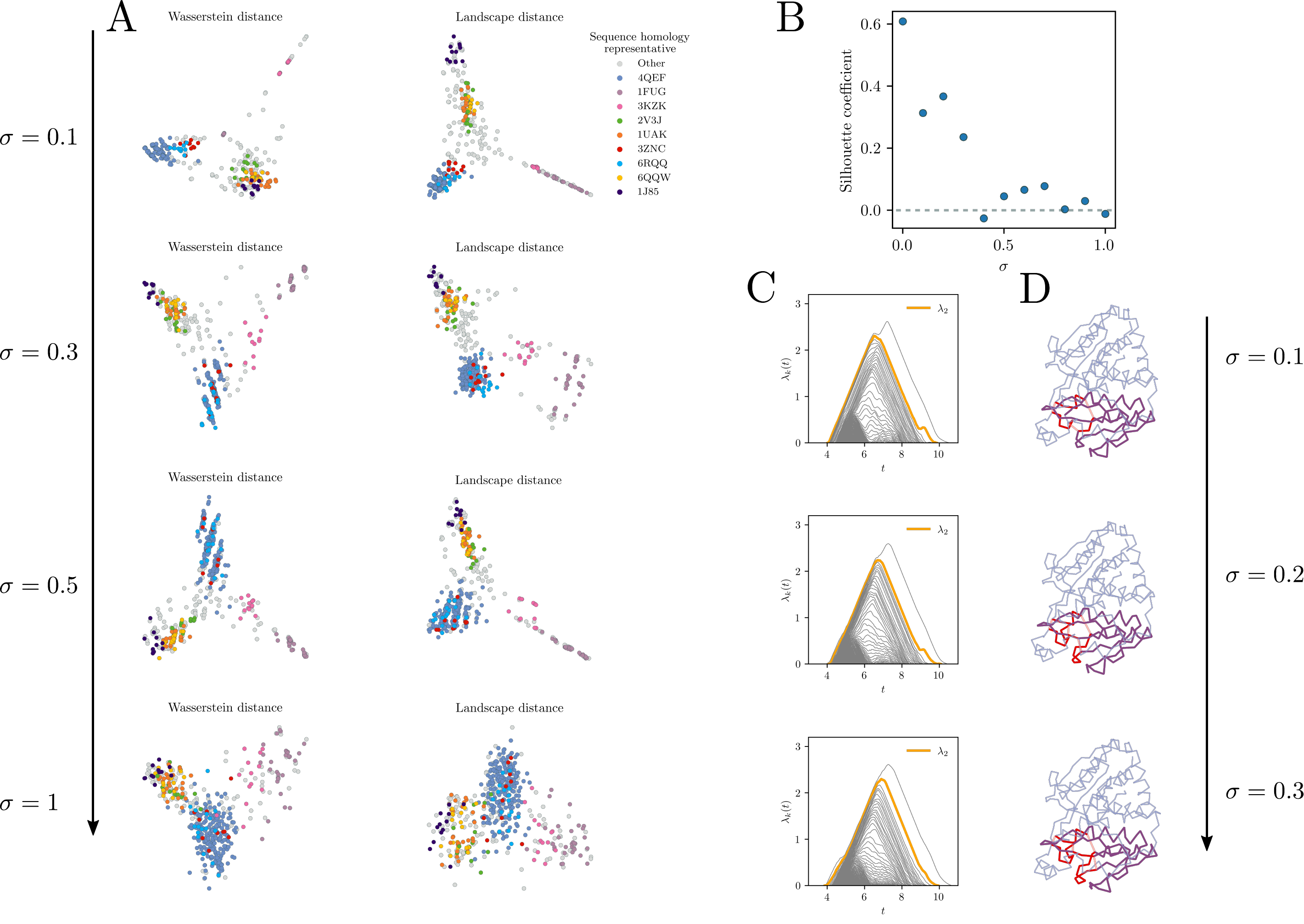}
\caption{\textbf{Sparse and noisy
dataset} \textbf{(A)} Evolution of the Isomap embeddings of the space of trefoil-knotted proteins equipped with the Wasserstein and landscape distances respectively as the standard deviation $\sigma$ of the Gaussian noise added increases. We see that the clustering separation holds, although it becomes less clear as the noise increases, and disappears for
$\sigma =1$. We also recover a clustering by depth type. The same plots with the depth-type colour code are available on the GitHub repository.\cite{ph-knotted-proteins}  \textbf{(B)} The silhouette coefficient\cite{rousseeuw1987silhouettes} of the sequence homology clustering with respect to the Isomap embedding on landscape distances as noise increases. \textbf{(C)} Evolution of the average landscape of the knotted AOTCases as the standard deviation $\sigma$ of the Gaussian noise added increases. We see that the $\lambda_2$ feature differentiating it from its unknotted counter-part stays for low values of $\sigma$ and is still almost visible for $\sigma = 0.3$. \textbf{(D)} The evolution of the cycle representing $c$ in the noisy 3KZK point cloud corresponding to the $\lambda_2$ peak. Note that in each of the cases the cycle is still positioned close to the knot core, and more specifically, close to the crossings responsible for the non-trivial entanglement. Furthermore, it still intersects the arc that needs to be pushed to untangle the curve.}
\label{fig2}
\end{figure}

\begin{table}[]
    \centering
    \include{pvalues.tex}
    \caption{The $p$-values output by the randomisation tests in SI \ref{sec:landscapecomps} when comparing the average persistence landscapes of each pair of sequence homology classes of trefoil-knotted protein chains.}
    \label{tab:pvalues}
\end{table}
\fi 
\end{appendices}


\begin{thebibliography}{10}

\bibitem{barbensi2021depth}
Agnese Barbensi and Dimos Goundaroulis.
\newblock $f$-distance of knotoids and protein structure.
\newblock {\em Proceedings of the Royal Society A}, 477:20200898, 02 2021.

\bibitem{barbensi2021topological}
Agnese Barbensi, Naya Yerolemou, Oliver Vipond, Barbara~I. Mahler, Pawel
  Dabrowski-Tumanski, and Dimos Goundaroulis.
\newblock A topological selection of folding pathways from native states of
  knotted proteins.
\newblock {\em Symmetry}, 13:1670, 2021.

\bibitem{bauer2021ripser}
Ulrich Bauer.
\newblock Ripser: efficient computation of {Vietoris-Rips} persistence
  barcodes.
\newblock {\em Journal of Applied and Computational Topology}, 5:391--423,
  2021.

\bibitem{ph-knotted-proteins}
Katherine Benjamin.
\newblock Code and data for {H}omology of {H}omologous {K}notted {P}roteins.
\newblock \url{https://github.com/katherine-benjamin/ph-knotted-proteins},
  2022.

\bibitem{berman2007worldwide}
Helen Berman, Kim Henrick, Haruki Nakamura, and John~L Markley.
\newblock The worldwide {P}rotein {D}ata {B}ank (wwpdb): ensuring a single,
  uniform archive of pdb data.
\newblock {\em Nucleic Acids Research}, 35(suppl\_1):D301--D303, 2007.

\bibitem{bubenik2015statistical}
Peter Bubenik.
\newblock Statistical topological data analysis using persistence landscapes.
\newblock {\em Journal of Machine Learning Research}, 16:77--102, 2015.

\bibitem{bubenik2017persistence}
Peter Bubenik and D\l{}otko Pawe\l{}.
\newblock A persistence landscapes toolbox for topological statistics.
\newblock {\em Journal of Symbolic Computation}, 78:91--114, 2017.

\bibitem{capraro2016untangling}
Dominique~T Capraro and Patricia~A Jennings.
\newblock Untangling the influence of a protein knot on folding.
\newblock {\em Biophysical Journal}, 110(5):1044--1051, 2016.

\bibitem{celoria2021statistical}
Daniele Celoria and Barbara~I. Mahler.
\newblock A statistical approach to knot confinement via persistent homology.
\newblock {\em arXiv preprint arXiv:2108.03034}, 2021.

\bibitem{cohen2007stability}
David Cohen-Steiner, Herbert Edelsbrunner, and John Harer.
\newblock Stability of persistence diagrams.
\newblock {\em Discrete \& Computational Geometry}, 37(1):103--120, 2007.

\bibitem{dabrowski2019knotprot}
Pawel Dabrowski-Tumanski, Pawel Rubach, Dimos Goundaroulis, Julien Dorier,
  Piotr Sulkowski, Kenneth~C Millett, Eric~J Rawdon, Andrzej Stasiak, and
  Joanna~I Sulkowska.
\newblock {KnotProt 2.0: a database of proteins with knots and other entangled
  structures}.
\newblock {\em Nucleic Acids Research}, 47(D1):D367--D375, 2019.

\bibitem{dabrowski2016}
Pawel Dabrowski-Tumanski, Andrzej Stasiak, and Joanna~I Sulkowska.
\newblock In search of functional advantages of knots in proteins.
\newblock {\em PloS one}, 11(11):e0165986--e0165986, 11 2016.

\bibitem{dorier2018knoto}
Julien Dorier, Dimos Goundaroulis, Fabrizio Benedetti, and Andrzej Stasiak.
\newblock Knoto-id: a tool to study the entanglement of open protein chains
  using the concept of knotoids.
\newblock {\em Bioinformatics}, 34(19):3402--3404, 2018.

\bibitem{dorier2021open}
Julien Dorier, Dimos Goundaroulis, Eric~J Rawdon, and Andrzej Stasiak.
\newblock Open knots.
\newblock {\em Encyclopedia of knot theory}, 2021.

\bibitem{edelsbrunner2010computational}
Herbert Edelsbrunner and John Harer.
\newblock {\em Computational topology: an introduction}.
\newblock American Mathematical Soc., 2010.

\bibitem{gameiro2014topological}
Marcio Gameiro, Yasuaki Hiraoka, Shunsuke Izumi, Miroslav Kramár, Konstantin
  Mischaikow, and Vidit Nanda.
\newblock A topological measurement of protein compressibility.
\newblock {\em Japan Journal of Industrial and Applied Mathematics}, 32:1--17,
  03 2014.

\bibitem{ghrist2008barcodes}
Robert Ghrist.
\newblock Barcodes: the persistent topology of data.
\newblock {\em Bulletin of the American Mathematical Society}, 45(1):61--75,
  2008.

\bibitem{goundaroulis2020chromatin}
Dimos Goundaroulis, Erez~Lieberman Aiden, and Andrzej Stasiak.
\newblock Chromatin is frequently unknotted at the megabase scale.
\newblock {\em Biophysical Journal}, 118(9):2268--2279, 2020.

\bibitem{goundaroulis2020knotoids}
Dimos Goundaroulis, Julien Dorier, and Andrzej Stasiak.
\newblock Knotoids and protein structure.
\newblock {\em Topology and Geometry of Biopolymers}, 746:185, 2020.

\bibitem{hamilton2021persistent}
Wesley Hamilton, J.~E. Borgert, T.~Hamelryck, and J.S. Marron.
\newblock Persistent topology of protein space.
\newblock {\em arXiv preprint arXiv:2102.06768}, 2021.

\bibitem{hatcher}
Allen Hatcher.
\newblock {\em Algebraic Topology}.
\newblock Cambridge University Press, 2002.

\bibitem{henselmanghristl6}
G.~{Henselman} and R.~{Ghrist}.
\newblock {Matroid Filtrations and Computational Persistent Homology}.
\newblock {\em arXiv preprint arXiv:1606.00199}, June 2016.

\bibitem{hiraoka2016hierarchical}
Yasuaki Hiraoka, Takenobu Nakamura, Akihiko Hirata, Emerson~G Escolar, Kaname
  Matsue, and Yasumasa Nishiura.
\newblock Hierarchical structures of amorphous solids characterized by
  persistent homology.
\newblock {\em Proceedings of the National Academy of Sciences},
  113(26):7035--7040, 2016.

\bibitem{jackson2020there}
Sophie~E Jackson.
\newblock Why are there knots in proteins?
\newblock {\em Topology and Geometry of Biopolymers}, 746:129, 2020.

\bibitem{jackson2017fold}
Sophie~E Jackson, Antonio Suma, and Cristian Micheletti.
\newblock How to fold intricately: using theory and experiments to unravel the
  properties of knotted proteins.
\newblock {\em Current Opinion in Structural Biology}, 42:6--14, 2017.

\bibitem{king2007identification}
Neil~P King, Eric~O Yeates, and Todd~O Yeates.
\newblock Identification of rare slipknots in proteins and their implications
  for stability and folding.
\newblock {\em Journal of Molecular Biology}, 373(1):153--166, 2007.

\bibitem{kovacev2016using}
Violeta Kovacev-Nikolic, Peter Bubenik, Dragan Nikoli\'{c}, and Giseon Heo.
\newblock Using persistent homology and dynamical distances to analyze protein
  binding.
\newblock {\em Statistical Applications in Genetics and Molecular Biology},
  15(1):19--38, 2016.

\bibitem{mallam2010experimental}
Anna~L Mallam, Joseph~M Rogers, and Sophie~E Jackson.
\newblock Experimental detection of knotted conformations in denatured
  proteins.
\newblock {\em Proceedings of the National Academy of Sciences},
  107(18):8189--8194, 2010.

\bibitem{mcguirl2020topological}
Melissa~R McGuirl, Alexandria Volkening, and Bj{\"o}rn Sandstede.
\newblock Topological data analysis of zebrafish patterns.
\newblock {\em Proceedings of the National Academy of Sciences},
  117(10):5113--5124, 2020.

\bibitem{millett2013identifying}
Kenneth~C Millett, Eric~J Rawdon, Andrzej Stasiak, and Joanna~I Su{\l}kowska.
\newblock Identifying knots in proteins.
\newblock {\em Biochemical Society Transactions}, 41(2):533--537, 2013.

\bibitem{roadmap}
Nina Otter, Mason~A Porter, Ulrike Tillmann, Peter Grindrod, and Heather~A
  Harrington.
\newblock A roadmap for the computation of persistent homology.
\newblock {\em EPJ Data Science}, 6:1--38, 2017.

\bibitem{scikit-learn}
F.~Pedregosa, G.~Varoquaux, A.~Gramfort, V.~Michel, B.~Thirion, O.~Grisel,
  M.~Blondel, P.~Prettenhofer, R.~Weiss, V.~Dubourg, J.~Vanderplas, A.~Passos,
  D.~Cournapeau, M.~Brucher, M.~Perrot, and E.~Duchesnay.
\newblock Scikit-learn: Machine learning in {P}ython.
\newblock {\em Journal of Machine Learning Research}, 12:2825--2830, 2011.

\bibitem{piejko2020folding}
Maciej Piejko, Szymon Niewieczerzal, and Joanna~I Sulkowska.
\newblock The folding of knotted proteins: Distinguishing the distinct behavior
  of shallow and deep knots.
\newblock {\em Israel Journal of Chemistry}, 60(7):713--724, 2020.

\bibitem{potestio2010knotted}
Raffaello Potestio, Cristian Micheletti, and Henri Orland.
\newblock Knotted vs. unknotted proteins: Evidence of knot-promoting loops.
\newblock {\em {PLOS} Computational Biology}, 6(7):1--10, 07 2010.

\bibitem{rolfsen2003knots}
Dale Rolfsen.
\newblock {\em Knots and links}, volume 346.
\newblock American Mathematical Society, 2003.

\bibitem{rose2021}
Yana Rose, Jose~M Duarte, Robert Lowe, Joan Segura, Chunxiao Bi, Charmi
  Bhikadiya, Li~Chen, Alexander~S Rose, Sebastian Bittrich, Stephen~K Burley,
  et~al.
\newblock RCSB Protein Data Bank: Architectural advances towards integrated
  searching and efficient access to macromolecular structure data from the PDB
  archive.
\newblock {\em Journal of Molecular Biology}, 433(11):166704, 2021.

\bibitem{rousseeuw1987silhouettes}
Peter~J. Rousseeuw.
\newblock Silhouettes: A graphical aid to the interpretation and validation of
  cluster analysis.
\newblock {\em Journal of Computational and Applied Mathematics}, 20:53--65,
  1987.

\bibitem{siebert2017there}
Jonathan~T. Siebert, Alexey~N. Kivel, Liam~P. Atkinson, Tim~J. Stevens,
  Ernest~D. Laue, and Peter Virnau.
\newblock Are there knots in chromosomes?
\newblock {\em Polymers}, 9(8):317, 2017.

\bibitem{stolz2020geometric}
Bernadette~J Stolz, Jared Tanner, Heather~A Harrington, and Vidit Nanda.
\newblock Geometric anomaly detection in data.
\newblock {\em Proceedings of the National Academy of Sciences},
  117(33):19664--19669, 2020.

\bibitem{sulkowska2012conservation}
Joanna~I Su{\l}kowska, Eric~J Rawdon, Kenneth~C Millett, Jose~N Onuchic, and
  Andrzej Stasiak.
\newblock Conservation of complex knotting and slipknotting patterns in
  proteins.
\newblock {\em Proceedings of the National Academy of Sciences},
  109(26):E1715--E1723, 2012.

\bibitem{sulkowska2008stabilizing}
Joanna~I Su{\l}kowska, Piotr Su{\l}kowski, P~Szymczak, and Marek Cieplak.
\newblock Stabilizing effect of knots on proteins.
\newblock {\em Proceedings of the National Academy of Sciences},
  105(50):19714--19719, 2008.

\bibitem{sulkowska2020folding}
Joanna~Ida Sulkowska.
\newblock On folding of entangled proteins: knots, lassos, links and
  $\theta$-curves.
\newblock {\em Current Opinion in Structural Biology}, 60:131--141, 2020.

\bibitem{tenenbaum2000}
Joshua~B. Tenenbaum, Vin De~Silva, and John~C. Langford.
\newblock A global geometric framework for nonlinear dimensionality reduction.
\newblock {\em Science}, 290(5500):2319--2323, 2000.

\bibitem{gudhi}
{The GUDHI Project}.
\newblock {\em GUDHI User and Reference Manual}.
\newblock GUDHI Editorial Board, 3.4.1 edition, 2021.

\bibitem{tubiana2011probing}
Luca Tubiana, Enzo Orlandini, and Cristian Micheletti.
\newblock Probing the entanglement and locating knots in ring polymers: a
  comparative study of different arc closure schemes.
\newblock {\em Progress of Theoretical Physics Supplement}, 191:192--204, 2011.

\bibitem{vipond2021multiparameter}
Oliver Vipond, Joshua~A Bull, Philip~S Macklin, Ulrike Tillmann, Christopher~W
  Pugh, Helen~M Byrne, and Heather~A Harrington.
\newblock Multiparameter persistent homology landscapes identify immune cell
  spatial patterns in tumors.
\newblock {\em Proceedings of the National Academy of Sciences}, 118(41), 2021.

\bibitem{2020SciPy-NMeth}
Pauli Virtanen, Ralf Gommers, Travis~E. Oliphant, Matt Haberland, Tyler Reddy,
  David Cournapeau, Evgeni Burovski, Pearu Peterson, Warren Weckesser, Jonathan
  Bright, St{\'e}fan~J. {van der Walt}, Matthew Brett, Joshua Wilson, K.~Jarrod
  Millman, Nikolay Mayorov, Andrew R.~J. Nelson, Eric Jones, Robert Kern, Eric
  Larson, C~J Carey, {\.I}lhan Polat, Yu~Feng, Eric~W. Moore, Jake
  {VanderPlas}, Denis Laxalde, Josef Perktold, Robert Cimrman, Ian Henriksen,
  E.~A. Quintero, Charles~R. Harris, Anne~M. Archibald, Ant{\^o}nio~H. Ribeiro,
  Fabian Pedregosa, Paul {van Mulbregt}, and {SciPy 1.0 Contributors}.
\newblock {{SciPy} 1.0: Fundamental Algorithms for Scientific Computing in
  Python}.
\newblock {\em Nature Methods}, 17:261--272, 2020.

\bibitem{xia2014persistent}
Kelin Xia and Guo-Wei Wei.
\newblock Persistent homology analysis of protein structure, flexibility, and
  folding.
\newblock {\em International Journal for Numerical Methods in Biomedical
  Engineering}, 30(8):814--844, August 2014.

\bibitem{zomorodian2005computing}
Afra Zomorodian and Gunnar Carlsson.
\newblock Computing persistent homology.
\newblock {\em Discrete \& Computational Geometry}, 33(2):249--274, 2005.

\end{thebibliography}
\end{document}